\documentstyle[12pt]{article}

\def\be{\begin{equation}}
\def\ee{\end{equation}}
\def\bea{\begin{eqnarray}}
\def\eea{\end{eqnarray}}

\catcode`\@=11
\font\tenmsa=msam10
\font\sevenmsa=msam7
\font\fivemsa=msam5
\font\tenmsb=msbm10
\font\sevenmsb=msbm7
\font\fivemsb=msbm5
\newfam\msafam
\newfam\msbfam
\textfont\msafam=\tenmsa  \scriptfont\msafam=\sevenmsa
  \scriptscriptfont\msafam=\fivemsa
\textfont\msbfam=\tenmsb  \scriptfont\msbfam=\sevenmsb
  \scriptscriptfont\msbfam=\fivemsb

\def\hexnumber@#1{\ifnum#1<10 \number#1\else
 \ifnum#1=10 A\else\ifnum#1=11 B\else\ifnum#1=12 C\else
 \ifnum#1=13 D\else\ifnum#1=14 E\else\ifnum#1=15
F\fi\fi\fi\fi\fi\fi\fi}

\def\msa@{\hexnumber@\msafam}
\def\msb@{\hexnumber@\msbfam}
\mathchardef\blacktriangleright="3\msa@49
\mathchardef\blacktriangleleft="3\msa@4A
\catcode`\@=12

\def\bicross{\triangleright\!\!\!\blacktriangleleft}

\def\wpq{\widetilde W_+^z}
\def\w{\widetilde}
\def\wz{\tilde z}
\def\p{\varepsilon}
\def\1{\'{\i}}

\begin{document}

\thispagestyle{empty}

\hfill \

\
\vspace{2cm}

\begin{center} {\LARGE{\bf{Bicrossproduct structure of the}}}

 {\LARGE{\bf{null-plane quantum Poincar\'e algebra}}}
 \end{center}

\bigskip\bigskip\bigskip

\begin{center}
Oscar Arratia $^{1}$, Francisco J. Herranz $^{2}$ and
Mariano A. del Olmo $^{3}$
\end{center}

\begin{center}
{\it $^{1}$ Departamento de Matem\'atica Aplicada a la
Ingenier\'{\i}a,  \\
Universidad de Valladolid. E-47011, Valladolid, Spain\\
E. mail: oscarr$@$wmatem.eis.uva.es }
\end{center}

\begin{center}
{\it $^{2}$ Departamento de F\1sica, E.U. Polit\'ecnica,\\
Universidad de Burgos.   E-09006, Burgos, Spain\\
E. mail: fteorica$@$cpd.uva.es}
\end{center}

\begin{center}
{\it $^{3}$ Departamento de F\'{\i}sica Te\'orica,\\
Universidad de Valladolid. E-47011, Valladolid, Spain\\
E. mail: olmo$@$cpd.uva.es}
\end{center}

\bigskip\bigskip\bigskip

\begin{abstract}
A nonlinear change of basis allows to show that  the non-standard
quantum deformation of the (3+1) Poincar\'e algebra  has a
bicrossproduct structure. Quantum universal $R$-matrix, Pauli--Lubanski
and mass operators are presented in the new basis.
\end{abstract}


\newpage

{\bf 1.-} The aim of this letter is to prove that the non-standard
quantum deformation of the (3+1) Poincar\'e algebra \cite{null},
the so-called null-plane quantum Poincar\'e algebra, can be
endowed with a structure of bicrossproduct \cite{majid}. After the
proofs by Majid and Ruegg \cite{kappa} that the (3+1) $\kappa$-Poincar\'e
algebra \cite{Luka,Gill,Lukb} has a bicrossproduct structure,
and more recently by Azc\'arraga {\em et al} \cite{ck} that the
$q$-Poincar\'e in any dimension \cite{BHOS} has also this
kind of structure, it only remains to study if the same
bicrossproduct structure is exhibited by the (3+1) null-plane
quantum Poincar\'e. In \cite{AP} it has been showed
that the (1+1) null-plane quantum Poincar\'e \cite{BHOPS} also shares
this structure, however, this lower dimensional  case does not indicate
the procedure for the (3+1) case, i.e., the nonlinear change of basis
that allows to display the bicrossproduct structure. It is worthy to
note that in all the three mentioned deformations the formal
decomposition is the same, i.e.,
$$
U_q({\cal P}(3+1)) = U(so(3,1))^\beta\bicross_\alpha
U_q({\cal T}_4),
$$
following the same pattern of the classical  algebra or group
counterpart
$$
 P(3+1) = SO(3,1)\odot T_4 ,
$$
and with the sector of the translations  deformed (differently in
each case, of course) and the Lorentz transformation sector
non-deformed.


{\bf 2.-} The generators of the (3+1) Poincar\'e algebra ${\cal
P}(3+1)$ in the  so-called null-plane basis \cite{LS} are
\be
\{P_+,P_-,P_i,E_i,F_i,K_3,J_3; \ i=1,2\},\label{basisold}
\ee
where $P_+$, $P_-$, $E_i$ and
$F_i$   are expressed in  terms of the usual  kinematical
ones $\{H,P_l,K_l,J_l;\ l=1,2,3\}$ by
\be
\begin{array}{lll}
 {P_+}= (H+P_3)/2,  &\quad  {P_-}=H-P_3,  &\quad
 {E_1}= (K_1+J_2)/2 , \cr
  {F_1}=K_1-J_2,  &\quad  {F_2}=K_2+J_1,  &\quad   {E_2}=
(K_2-J_1)/2 .
\end{array}
\label{aa}
\ee
Hence, the Lie brackets of ${\cal P}(3+1)$ are
(hereafter $i,j=1,2$):
\be
\begin{array}{lll}
   [K_3,E_i]=E_i, &\qquad
[K_3,F_i]=-F_i, &\qquad [K_3,J_3]=0,\cr
 [J_3,E_i]=-\p_{ij3}E_j,
&\qquad [J_3,F_i]=-\p_{ij3}F_j, &\qquad [E_1,E_2]=0,\cr
 \multicolumn{3}{l} { [E_i,F_j]=\delta_{ij}K_3
+\p_{ij3}J_3,\qquad\quad
 [F_1,F_2]=0,}
\end{array}
\label{ab}
\ee
\be
[P_\mu,P_\nu]=0,\qquad \mu,\nu=+,-,1,2,
\label{ac}
\ee
\be
\begin{array}{lll}
 [K_3,P_+]=P_+, &\qquad [K_3,P_-]=-P_-, &\qquad
 [K_3,P_i]=0 ,\cr
 [J_3,P_i]=-\p_{ij3}P_j,
&\qquad    [J_3,P_+]=0,&\qquad [J_3,P_-]=0,\cr
[E_i,P_j]=\delta_{ij}P_+,
 & \qquad [E_i,P_-]=P_i, & \qquad [E_i,P_+]=0,\cr
  [F_i,P_j]=\delta_{ij}P_-,
&\qquad [F_i,P_+]= P_i, &\qquad [F_i,P_-]= 0,
\end{array}
\label{ad}
\ee
 where $\p_{ijk}$ is the completely skewsymmetric tensor.

The semidirect product structure of the (3+1)  Poincar\'e group,
isomorphic to $ISO(3,1)$, can be clearly pointed out. The six
generators $\{E_i,F_i,K_3,J_3\}$ close the  Lorentz subgroup
$SO(3,1)$  (\ref{ab}), while the four remaining
$\{P_+,P_-,P_i\}$ generate the abelian subgroup $T_4$
(\ref{ac}). Therefore, as it is well known,
$ISO(3,1)=SO(3,1)\odot T_4$.


{\bf 3.-} A triangular or non-standard quantum deformation of  ${\cal
P}(3+1)$ was introduced in \cite{null} in the null-plane framework above
mentioned, whose Hopf structure we rewrite here for sake of
completeness and to clarify our main result. Let us denote the
null-plane generators $X$ displayed in (\ref{basisold}), by  $\w
X$, and by $\tilde z$ the deformation parameter.

\noindent
{\sl Coproduct}:
\bea
&&\Delta(\w  X)=1\otimes \w  X+\w  X\otimes 1, \qquad
\mbox{for}\quad \w
X\in\{\w P_+,\w E_i,\w J_3\},\cr
&&\Delta(\w Y)=e^{-\wz\w P_+}\otimes
\w Y+\w Y\otimes e^{\wz\w P_+},\quad
\mbox{for}\quad \w Y\in\{\w P_-,\w P_i\},\cr
&& \Delta(\w F_1)=e^{-\wz\w P_+}\otimes \w F_1
+\w F_1\otimes e^{\wz\w P_+}
+\wz e^{-\wz\w P_+} \w E_1
\otimes \w P_- - \wz\w P_-\otimes \w E_1
e^{\wz\w P_+}\cr
 &&\qquad\qquad  +\wz e^{-\wz\w P_+} \w J_3
\otimes \w P_2 - \wz \w P_2\otimes \w J_3
e^{\wz\w P_+},\label{co3}\\ &&
\Delta(\w F_2)=e^{-\wz\w P_+}\otimes
\w F_2+\w F_2\otimes e^{\wz\w P_+}
+\wz e^{-\wz\w P_+} \w E_2
\otimes \w P_- - \wz\w P_-\otimes \w E_2
e^{\wz\w P_+}\cr &&\qquad\qquad  -\wz e^{-\wz\w P_+} \w J_3
\otimes \w P_1 + \wz \w P_1\otimes \w J_3
e^{\wz\w P_+},\cr && \Delta(\w K_3)=e^{-\wz\w P_+}\otimes
\w K_3+\w K_3\otimes e^{\wz\w P_+}
+\wz e^{-\wz\w P_+} \w E_1
\otimes \w P_1 - \wz\w P_1\otimes \w E_1
e^{\wz\w P_+}\cr &&\qquad\qquad
+\wz e^{-\wz\w P_+} \w E_2
\otimes \w P_2 - \wz \w P_2\otimes\w E_2
e^{\wz\w P_+}; \nonumber
\eea
{\sl Counit and Antipode}:
\be
\epsilon(\w X) =0;\quad
\gamma(\w X)=-e^{{3\wz}\w P_+}\ \w X\ e^{- {3\wz}\w P_+},\quad
\mbox{for $\w X\in
\{\w P_\pm,\w P_i,\w E_i,\w F_i,
\w K_3,\w J_3\}$};
\label{anti3}
\ee
{\sl Non-vanishing Lie brackets}:
\bea
&& [\w K_3,\w P_+]=\frac {\sinh \wz\w P_+}{\wz},\quad
[\w K_3,\w P_-]=-\w P_-
\cosh \wz\w P_+,\quad [\w K_3,\w E_i]=\w E_i\cosh \wz\w P_+,\cr
&& [\w K_3,\w F_1]=- \w F_1 \cosh \wz\w P_+  +  \wz
\w E_1 \w P_-\sinh \wz\w P_+ - \wz^2\w P_2\,\wpq,\cr
&& [\w K_3,\w F_2]=- \w F_2 \cosh \wz\w P_+
+  \wz  \w E_2 \w P_- \sinh \wz\w P_+  + \wz^2\w P_1\,\wpq,\cr
&& [\w J_3,\w P_i]=-\p_{ij3}\w P_j, \qquad  [\w J_3,\w E_i]=-\p_{ij3}\w
E_j,
\qquad [\w J_3,\w F_i]=-\p_{ij3}\w F_j,\cr
&& [\w E_i,\w P_j]=\delta_{ij}\frac {\sinh \wz\w P_+}{\wz}, \qquad
[\w F_i,\w P_j]=\delta_{ij}\w P_-\cosh\w  z\w P_+,\label{brack3}\\
&& [\w E_i,\w F_j]=\delta_{ij}\w K_3 +\p_{ij3}\w J_3\cosh \wz\w
P_+,\qquad [\w P_+,\w F_i]=-\w P_i,\cr && [\w F_1,\w F_2]=\wz^2
\w P_- \wpq +
\wz  \w P_- \w J_3 {\sinh \wz\w P_+}, \qquad [\w P_-,\w E_i]=-\w
P_i,
\nonumber
\eea
where   $\wpq$ is  a component of the deformed
Pauli--Lubanski vector defined as
\be
\wpq= \w E_1\w P_2-\w E_2\w P_1 +
\w J_3\frac{\sinh \wz\w P_+}{\wz}.
\ee


{\bf 4.-} In the sequel we
show that this quantum algebra has a bicrossproduct structure
\cite{majid}. Let us consider the map  defined by:
\bea
&&P_+={\w  P}_+, \qquad E_i={\w  E}_i, \qquad
J_3={\w  J}_3,\qquad  z=2\tilde z, \cr
&&P_-=e^{-\tilde z{\w  P}_+}{\w  P}_-, \qquad
P_i=e^{-\tilde z{\w  P}_+}{\w  P}_i, \cr
&&F_1=e^{-\tilde z{\w  P}_+}({\w  F}_1-
\tilde z {\w  E}_1{\w  P}_- -
\tilde z {\w  J}_3{\w  P}_2) \label{ba},\\
&&F_2=e^{-\tilde z{\w  P}_+}({\w  F}_2-
\tilde z {\w  E}_2{\w  P}_- +
\tilde z {\w  J}_3{\w  P}_1 ), \cr
&&K_3=e^{-\tilde z{\w  P}_+}({\w  K}_3-
\tilde z {\w  E}_1{\w  P}_1 -
\tilde z {\w  E}_2{\w  P}_2 ).
\nonumber
\eea
By applying (\ref{ba}) to
the Hopf algebra $U_{\wz}({\cal P}(3+1))$, whose relations appear
displayed in expressions (\ref{co3})--(\ref{brack3}), we get  the
Hopf algebra $U_z({\cal P}(3+1))$, characterized by the  following
coproduct, counit, antipode and commutation relations:
\bea
&&\Delta(X)=1\otimes X+X\otimes 1, \quad
X\in\{P_+,E_i,J_3\}, \cr
&&\Delta(Y)=e^{-zP_+}\otimes Y+Y\otimes 1,\qquad
 Y\in\{P_-,P_i\}, \cr
&& \Delta(F_1)=e^{-zP_+} \otimes F_1
+F_1\otimes 1 -  zP_-\otimes E_1
 -  z P_2\otimes J_3, \label{bb}\\
&& \Delta(F_2)=e^{-zP_+}\otimes F_2+F_2\otimes 1
 -  zP_-\otimes E_2   +
 z P_1\otimes J_3, \cr
&& \Delta(K_3)=e^{-zP_+}\otimes K_3+K_3\otimes 1
 -  zP_1\otimes E_1   -
 z P_2\otimes E_2;
\nonumber
\eea
\be
\epsilon(X) =0, \qquad   X\in \{P_\pm,P_i,E_i,F_i,K_3,J_3\};
\label{bc}
\ee
\bea
&&\gamma(X)=-X,\qquad  X\in \{P_+,E_i,J_3\}, \cr
&&\gamma(Y)=-e^{zP_+}Y, \qquad Y\in \{P_-,P_i\}, \cr
&&\gamma(F_1)=-e^{zP_+}(F_1+zP_- E_1 + z P_2 J_3), \label{bd}\\
&&\gamma(F_2)=-e^{zP_+}(F_2+zP_- E_2 - z P_1 J_3),\cr
&&\gamma(K_3)=-e^{zP_+}(K_3+zP_1 E_1 +  z P_2 E_2);
\nonumber\eea
\be
\begin{array}{lll}
   [K_3,E_i]=E_i, &\qquad
[K_3,F_i]=-F_i, &\qquad [K_3,J_3]=0,\cr
 [J_3,E_i]=-\p_{ij3}E_j,
&\qquad [J_3,F_i]=-\p_{ij3}F_j, &\qquad [E_1,E_2]=0,\cr
 \multicolumn{3}{l}{ [E_i,F_j]=\delta_{ij}K_3 +\p_{ij3}J_3,
\qquad\quad
 [F_1,F_2]=0,}
\end{array}
\label{be}
\ee
\be
[P_\mu,P_\nu]=0,\qquad \mu,\nu=+,-,1,2,
\label{bf}
\ee
\be
\begin{array}{ll}
\displaystyle
 [K_3,P_+]=\frac {1-e^{-zP_+}}{z}, &\qquad\displaystyle [K_3,P_-]=-P_-
-\frac z2(P_1^2+P_2^2),\cr
 [K_3,P_i]=(e^{-zP_+}-1)P_i, & \qquad [J_3,P_+]=0 ,
 \qquad    [J_3,P_-]=0 ,\cr
  [J_3,P_i]=-\p_{ij3}P_j,  &\qquad [E_i,P_-]=P_i,
\qquad [E_i,P_+]=0,\cr
\displaystyle [E_i,P_j]=\delta_{ij}\frac {1-e^{-zP_+}}{z},&\qquad
[F_i,P_+]= P_i,  \qquad [F_i,P_-]= -zP_iP_-, \cr
 \multicolumn{2}{l}{
\displaystyle [F_i,P_j]=-zP_iP_j+\delta_{ij}(e^{-zP_+} P_- +\frac
z2(P_1^2+P_2^2)) .}
\label{bg}
\end{array}
\ee

Note that the translation generators $\{P_+,P_-,P_i\}$ define a
commutative but non-cocommu\-tative Hopf subalgebra  of $U_z({\cal
P}(3+1))$ denoted $U_z({\cal T}_4)$, and the Lorentz sector is
non-deformed at the algebra level.

{\bf 5.-} Let us consider now the non-deformed Lorentz Hopf
algebra $U(so(3,1))$ spanned by the  generators $\{E_i,F_i,K_3,J_3\}$
with classical commutation rules
(\ref{ab}) and primitive coproduct:  $\Delta(X)=1\otimes X+X\otimes
1$. We   define a right action
\be
\alpha: U_z({\cal T}_4)\otimes U(so(3,1))\to U_z({\cal T}_4)
\ee
 as
\be
\alpha(X\otimes Y)\equiv X\lhd Y:=[X,Y], \qquad X\in\{P_\pm,P_i\},
\quad Y\in \{E_i,F_i,K_3,J_3\};
\label{bh}
\ee
explicitly
\be
\begin{array}{ll}
\displaystyle
 \alpha (P_+\otimes K_3) =\frac {e^{-zP_+}-1}{z},
&\qquad\displaystyle
\alpha (P_-\otimes  K_3) = P_-  +\frac z2(P_1^2+P_2^2),\cr
\alpha (P_i\otimes  K_3)  =(1-e^{-zP_+} )P_i, & \qquad
\alpha (P_+\otimes  J_3) =0,
 \qquad \quad \alpha (P_-\otimes  J_3)=0, \cr
  \alpha (P_i\otimes  J_3) =\p_{ij3}P_j,  &\qquad
 \alpha (P_-\otimes E_i) =-P_i,
\qquad \alpha ( P_+\otimes E_i) =0,\cr
\displaystyle \alpha (P_i\otimes E_j) =\delta_{ij}\frac {
e^{-zP_+}-1 }{z},&\qquad   \alpha (P_+\otimes F_i) =- P_i , \qquad
 \alpha (P_-\otimes F_i) = zP_iP_- ,\cr
 \multicolumn{2}{l}{
\displaystyle \alpha (P_i\otimes F_j)
 = zP_iP_j-\delta_{ij}(e^{-zP_+} P_-
+\frac z2(P_1^2+P_2^2)) .}
\label{bi}
\end{array}
\ee
Also we   define a left coaction
\be
\beta:   U(so(3,1))\to U_z({\cal
T}_4)\otimes U(so(3,1))
\ee
by
\bea
&&\beta(J_3)=1\otimes J_3,\qquad \beta(E_i)=1\otimes E_i, \cr
&& \beta(F_1)=e^{-zP_+} \otimes F_1
 -  zP_-\otimes E_1
 -  z P_2\otimes J_3, \label{ca}\\
&& \beta(F_2)=e^{-zP_+}\otimes F_2   -  zP_-\otimes E_2   +
 z P_1\otimes J_3, \cr
&& \beta(K_3)=e^{-zP_+}\otimes K_3
 -  zP_1\otimes E_1 -
 z P_2\otimes E_2   .
\nonumber
\eea

It can be shown that the right action $\alpha$ and left coaction
$\beta$ just introduced fulfill  the compatibility conditions
\cite{majid} in such manner $(U_z({\cal T}_4),\alpha)$ is a right
$U(so(3,1))$-module algebra and $(U(so(3,1)),\beta)$ is a left
$U_z({\cal T}_4)$-comodule coalgebra.
We  summarize the previous discussion in the following theorem,
which is the main result of this letter together with the nonlinear
basis change (\ref{ba}).

\noindent
{\bf Theorem.}
{\em The null-plane
quantum Poincar\'e algebra has the bicrossproduct structure}
\be
U_z({\cal P}(3+1)) = U(so(3,1))^\beta\bicross_\alpha
U_z({\cal T}_4) .
\label{cb}
\ee

{\bf 6.-} We would like to stress that the map (\ref{ba}) is
invertible, so:
\bea
&&{\w  P}_+=P_+, \qquad {\w  E}_i=E_i ,\qquad
{\w  J}_3=J_3 , \qquad \tilde z=z/2,\cr
&&{\w  P}_-=e^{ z P_+/2} P_-, \qquad
{\w  P}_i=e^{  z P _+/2} P_i,\cr
&&{\w   F}_1=e^{ z P_+/2} (F_1 +
  z (E_1P_- +  J_3 P_2)/2), \label{ma}\\
&&{\w   F}_2=e^{ z P_+/2} (F_2 +
  z (E_2P_- -   J_3 P_1)/2) ,\cr
&&{\w   K}_3=e^{ z P_+/2} (K_3 +
  z (E_1P_1 +  E_2 P_2)/2).
\nonumber
\eea
This fact  can be applied to reproduce in the bicrossproduct basis
the physically relevant operators introduced in \cite{null} such
as Casimirs, spin, Hamiltonians and position operators. In particular,
the deformed square of the mass $M_z^2$ is now
\be
M_z^2= 2P_-\frac{e^{zP_+}-1}z -
(P_1^2 + P_2^2 ) e^{zP_+} ,
\label{mb}
\ee
and the square of the
Pauli--Lubanski operator $W_z^2$ turns out to be
\be
W_z^2=(W^z_{13})^2 +(W^z_{23})^2 +
\cosh(zP_+/2)\left( W^z_{+}W^z_{-} + W^z_{-} W^z_{+}\right)
-z^2 M_z^2(W^z_{+})^2/4,
\label{mc}
\ee
where
\bea
&&W^z_{i3}=K_3P_ie^{ z P_+} + E_i P_- - F_i\frac{e^{ z P_+}-1}{z}
+ \frac z2 (E_1P_1+E_2P_2)P_i e^{ z P_+}\cr
&&\qquad\qquad +
(-1)^i J_3P_{3-i}\frac{e^{ z P_+}-1}{2},\quad\qquad i=1,2,\cr
&&W^z_{-}= (F_1P_2-F_2P_1)e^{ z P_+}+J_3P_-\frac{e^{ z P_+}+1}2
+\frac z2 (E_1P_2-E_2P_1) P_- e^{ z P_+}\cr
&&\qquad \qquad
+\frac z2 J_3(P_1^2+P_2^2)e^{ z P_+}, \cr
&&W^z_{+}= (E_1P_2-E_2P_1)e^{ z P_+/2} +
J_3\frac{\sinh(zP_+/2)}{z/2} .
\label{md}
\eea
The second order Casimir (\ref{mb}) would give rise to a deformed
Schr\"odinger equation in the same way as in \cite{null,ba}, while
the  Pauli--Lubanski vector (\ref{md}) would allow  to derive quantum
Hamiltonians and spin operators.  However, it is clear that
although the Hopf algebra structure of $U_z({\cal P}(3+1))$ is
rather simplified in this new basis,  the associated  operators
adopt a much more complicated form than the original ones.

The map (\ref{ba}) resembles the one given in
\cite{rmatrix} which allowed to deduce a (factorized) null-plane
quantum universal $R$-matrix:  both mappings are related by the
interchange    $e^{-\tilde z{\w  P}_+} \leftrightarrow
e^{\tilde z{\w  P}_+}$. Hence, the universal $R$-matrix reads now
\bea
&&{\cal R}=\exp\{z E_2\otimes e^{z P_+} P_2\}
\exp\{z E_1\otimes e^{z P_+} P_1\}
\exp\{- z P_+\otimes e^{z P_+} K_3\}\cr
&&\qquad \times \exp\{ z e^{z P_+} K_3\otimes P_+\}
\exp\{- z e^{z P_+} P_1\otimes E_1\}
\exp\{- z e^{z P_+} P_2\otimes  E_2\}.
\eea

Therefore, each basis seems
to be useful for a specific  purpose and we do not find any
privileged basis to express  the whole quantum  Poincar\'e algebra
together its associated elements (universal $R$-matrix, quantum
Casimirs, etc.).

Finally  to mention that the null-plane case in $(2+1)$ dimensions
\cite{ba} also exhibits  this bicrossproduct  structure. It looks
interesting to profit this bicrossproduct structure of the quantum
algebras in order to study  by duality the corresponding quantum
groups. Work in this direction is in progress and will be published
elsewhere.


\noindent
{\section*{Acknowledgments}}

This work has been partially supported by the
Direcci\'on General de Ense\~nanza Superior (DGES grant PB95--0719)
and by DGICYT (project PB94--1115) from the Ministerio de
Educaci\'on  y Cultura of Spain.

\bigskip



\end{document}